\documentclass[12pt]{article}
\usepackage{amsmath}
\usepackage{amssymb}
\usepackage{amsthm}
\usepackage{amsfonts}

\numberwithin{equation}{section}

\newcommand{\R}{{\mathbb R}}
\newcommand{\Z}{{\mathbb Z}}

\newcommand{\1}[1]{\mathbf 1{#1}}

\def\bR{\mathbf R} 
\def\bQ{\mathbf Q}\def\vphi{\varphi} \def\R{\mathbb R}
\def\rB{{\rm B}} \def\rT{{\rm T}}
\def\gam{\gamma} \def\olam{\overline\lambda} \def\omu{\overline\mu}
\def\ogam{\overline\gamma}

\def\uy{{\underline y}} \def\bbZ{{\mathbb Z}} \def\gam{{\gamma}}
\def\eps{\epsilon}

\def\Lam{{\Lambda}}
\def\lam{{\lambda}}

\def\unn{{\underline n}}
\def\ue{{\underline e}}

\def\diy{\displaystyle}
\def\beq{\begin{equation}}
\def\eeq{\end{equation}}

\newtheorem{theorem}{Theorem}[section]

\theoremstyle{definition}

\begin{document}\def\cl{\centerline}

\cl{\LARGE{\bf Random walks}}
\vskip .5 truecm
\cl{\LARGE{\bf in a queueing network environment}}
\vskip 1 truecm

\cl{\Large{\bf M Gannon$^1$, E Pechersky$^{1,2}$, Y Suhov$^{2-5}$, A Yambartsev$^{1,6}$}}

\date{\today}

\vskip .5 truecm

{\bf Abstract.} We propose a class of models of random walks in a random environment where
an exact solution can be given for a stationary distribution. The environment is cast in terms of 
a Jackson/Gordon-Newell network although alternative interpretations are possible. The main 
tool is the detailed balance equations. The difference with earlier works is that the 
position of the random walk influences the transition intensities of the network environment
and vice versa, creating strong correlations. The form of the stationary distribution is closely 
related to the well-known product-formula.

\vskip .5 truecm



\section{Introduction} This paper introduces exactly solvable reversible models of a random walk interacting 
with a random environment of  a queueing-network type.
The environment stems from a symmetric Jackson network or its closed version (a Gordon--Newell network); cf. \cite{J1,J2, GoNe}. The walking particles can be  interpreted as distinguished
customers (DCs). Depending on the site of the network a DC is in, he/she may like staying in 
the site if the queue size is large (or small) and in turn encourages more\\
\vskip .1 truecm
--------------------------------
\vskip .5 truecm
$^1$ IME, University of Sao Paulo (USP), Sao Paulo 05508-090,
SP, Brazil

$^2$ IITP, Moscow, 127994, RF

$^3$ ICMC, University of Sao Paulo (USP), Sao Carlos 13566-590 SP, Brazil

$^4$ Statistical Laboratory, DPMMS, University of Cambridge, Cambridge CB3 0WB, UK

$^5$ Math Department, Penn State University, University Park, State College, PA 16802, USA

$^6$ Department of Pharmacy, Oregon State University, Corvallis, OR 97331-4605, USA
\newpage

\noindent (or less) tasks to come to this site. Further development is where DCs interact with each other: 
here we consider
the case of a symmetric simple exclusion; the zero-range model is also included. By analyzing the 
detailed balance equations (DBEs),  the equilibrium 
distribution is derived, closely related to the product-form distribution; cf. \cite{BCMP, N, Se2}. 
Our approach can be considered as a 
further development of earlier papers \cite{Ec, Zhu} where possible 
applications have been outlined, in the area of communication networks.

Other areas of applications may cover problems of  random trapping/localization and condensation; cf.       \cite{CW, EH, GSS, Szn}, and the bibliography therein. 

An immediate problem would be to identify 
a Lyapunov function assessing the speed of convergence to the stationary (equilibrium) distribution
under sub-criticality conditions. Cf.  \cite{M, IT} and references therein.

A possible step would be to introduce similar models in quasi-reversible 
setting \cite{Ke}. In this direction, it would be interesting to treat a wider class of service rules 
for tasks and the DC. Viz., a processor-sharing discipline seems a natural choice, involving a 
branching formalism; see, e.g., \cite{Se1}.  

Our models also show some potential in the direction of Markov processes with local interaction. An 
immediate goal can be to develop models in an infinite-volume configuration space \cite{Spi1}, \cite{Do1}, \cite{Do2}. 

\vskip .5 truecm\def\obeta{\overline\beta}

\section{Description of a basic model}
\subsection{ A symmetric open Jackson network} 
As a model for an environment, we take a symmetric and homogeneous Jackson job-shop 
network; see \cite{J1,J2}. The model is defined by the following ingredients.

\begin{description}
\item[(a)]\label{1.0} A finite collection $\Lam$ of sites with a 
single-server system assigned to each $k\in\Lam$. 
\item[(b)] Two  positive numbers: 
$\lam>0$, the intensity of an exogenous input flow to a given site,  and $\mu >0$, 
the intensity of the flow from a  given site out of the network.
\item[(c)] A non-negative symmetric matrix of transmission intensities
$\;\rB=(\beta_{ik},\,i,k\in\Lam )$: 
\beq\label{2.1}\beta_{ik}=\beta_{ki},\;\;\;\beta_{ik}\geq 0,\;\;\;\beta_{ii}=0.\eeq
\end{description}
 
The value $\lam$ gives the rate of independent Poisson processes of exogenous tasks arriving 
at sites in $\Lam$ and
$\omu_i=\mu+\obeta_i$ is the rate of servicing the queue at site $i\in\Lam$ where ${\obeta}_i
=\sum\limits_{k\in\Lam}\beta_{ik}$. After completing service at site $i$,
a task leaves the network with probability $\mu/\omu_i$ and jumps to site $k$ with probability $\beta_{ik}/\omu_i$.
Condition $\beta_{ii}=0$ is used for convenience. In contrast, the symmetry property $\beta_{ik}=\beta_{ki}$ 
in \eqref{2.1} is essential but rather restrictive and hopefully could be weakened in future. 

The above description gives rise to a continuous-time Markov process (MP) with states 
$\unn=(n_i,\,i\in\Lam)\in\Z_+^\Lam$ where $n_i\in\Z_+:=\{0,1,2,\ldots\}$.  The generator 
matrix $\bQ=(Q(\unn ,\unn^\prime))$ of the process has non-zero entries corresponding to 
the following transitions:
\beq\label{2.2}\begin{array}{ll}
{Q}(\unn ,\unn +\ue^i)=\lam&\hbox{an exogenous arrival of a task at site $i$,}\\
{Q}(\unn ,\unn -\ue^i)=\mu{\mathbf 1}(n_i\geq 1)&\hbox{a task exits 
from site $i$ out of the network,}\\
{Q}(\unn ,\unn +\ue^{i\to k})=\beta_{ik}{\mathbf 1}(n_i\geq 1)&\hbox{a task jumps from site $i$ 
to $k$.}\end{array}
\eeq
Here $\ue^i=(e^i_l,\,l\in\Lam)\in\Z_+^\Lam$ \ has \ $e^i_l=\delta_{il}$, $\ue^{i\to k}$ denotes the
difference $\ue^k-\ue^i$, \ and  ${\mathbf 1}$ stands 
for the indicator.

Assuming the sub-criticality condition 
$\diy\frac\lambda\mu<1$, 
the invariant measure $\pi$ is the product of geometric distributions with parameter $\lam/\mu$.
Here the probability $\pi (\unn )$ of 
having $n_i$ tasks at sites $i\in\Lam$ is of the form  
\beq\label{2.3}\pi (\unn )=\left(1-\frac{\lam}{\mu}\right)^{\#\,\Lam}\prod\limits_{i\in\Lam}\left(\frac\lambda
\mu\right)^{n_i},\;\;\unn=(n_i)\in\Z_+^\Lam.\eeq
In fact, assuming that matrix $\rB$ is irreducible (i.e., $\rB^t$ has strictly positive entries
for some $t\in\bbZ_+$), the process is positive recurrent (PR).

Eqn \eqref{2.3} follows from the DBEs for probabilities
$\pi (\unn )$, $\unn\in\Z_+^\Lam$: 
$$\begin{array}{c}
\pi (\unn)Q(\unn,\unn+\ue^i)=\pi (\unn+\ue^i)Q(\unn+\ue^i,\unn),\;\;
\pi (\unn)Q(\unn,\unn-\ue^j)=\pi (\unn-\ue^j)Q(\unn-\ue^j,\unn),\;n_j\ge q,\\
\pi (\unn)Q(\unn,\unn+\ue^{k\to l})=\pi (\unn+\ue^{k\to l})Q(\unn+\ue^{k\to l},\unn),\;n_k\geq 1,
\end{array}$$
which are easy to check. Note that the probabilities $\pi (\unn )$ in \eqref{2.3} do not refer to matrix\ $\rB$.

In the whole paper, $\Lam$ stands for a finite set, and sites of $\Lam$ are marked by $i,j,k,l,p,q,r,s$ and 
$j^\prime$.

\subsection{Random walk in a Jackson environment} We now give the description of the model with interaction. A new ingredient is the presence of a random particle (a distinguished customer (DC)) walking over set
$\Lam$. Parameters of the Jackson network are changed only for a site where the DC is located;
we call it a \textit{loaded} site and denote by $j$. A site $i\neq j$ (free of a DC) is called {\it unloaded}. In general, $i$, $k$ 
and $l$ denote sites in $\Lam$ which may be loaded or unloaded.

In addition to  $\lam$, $\mu$ and $\rB =(\beta_{ik})$ (see items (b) and (c) in Section 2.1), we need  more 
ingredients and rules. The queue of tasks at the loaded site $j$ has an exogenous arrival intensity $e^\varphi\lam$ where
$\varphi\in \R$, while for the remaining sites the intensity remains $\lam$. The exit flow intensity from all
sites of the network equals $\mu$ as before. The intensities $\beta_{ik}$ for $i\neq j\neq k$ 
(task jumps from one unloaded site to another) are as before (and satisfy \eqref{2.1}). Next, we introduce 
symmetric matrices $\Theta =(\theta_{ik},\,i,k\in\Lam )$ and $\rT=(\tau_{ik},\;i,k\in\Lam )$ where
\beq\label{2.4}
\theta_{ik}=\theta_{ki},\;\;\theta_{ik}\geq 0,\;\;\theta_{ii}=0;\;\;\tau_{ik}=\tau_{ki},\;\;\tau_{ik}\geq 0,\;\;\tau_{ii}=0.\eeq
(Again, condition $\theta_{ii}=\tau_{ii}=0$ is used for convenience.)
For the loaded site $j$, the intensity 
of the task flow from $j$ to $k\neq j$ equals $e^\varphi\theta_{jk}$. For an unloaded site $k\neq j$, the intensity
of the task flow from $k$ to $j$ equals $\theta_{kj}$. Finally, the intensity of a leap 
of the DC from a loaded site $j$ where there are $n_j$ tasks to a site $j^\prime\neq j$ is taken to
 be $e^{-\vphi n_j}\tau_{jj^\prime}$. As with intensities $\beta_{ik}$, symmetry equations $\theta_{ik}=\theta_{ki}$
 and $\tau_{ik}=\tau_{ki}$ are essential (in modified forms, they reappear in the rest of the paper), and it would be 
 interesting to replace these conditions with weaker ones.
 
Matrices $\Theta$ and $\rT$ and parameter $\varphi\in\R$ are additional ingredients of 
the model. (In future we also refer to $\rB$, $\Theta$ and $\rT$ as arrays, or collections 
of intensities.) The state of the emerging MP is a pair $(j,\unn )$ where
$j\in\Lam$ and $\unn =(n_i,\,i\in\Lam)\in\Z_+^\Lam$. As was  said, the entry $j$ indicates the loaded site, and 
$n_i$, as before, gives the number of tasks at site $i\in\Lam$. In accordance with the above description, 
the generator  ${\bR}=\Big\{R\big[(j,\unn );(j^\prime ,\unn^\prime)\big]\Big\}$ 
has the following entries (here and below, zero-rate transitions are not shown):
\beq\label{2.6}\begin{array}{ll}
{R}[(j,\unn) ;(j,\unn +\ue^k)]=\lam e^{\varphi\delta_{jk}}&\mbox{exogenous arrival of a task,}\\
{R}[(j,\unn) ;(j,\unn -\ue^k)]=\mu{\mathbf 1}(n_k\ge 1)&\mbox{exit of a task out of the network,}\\
{R}[(j,\unn) ;(j,\unn +\ue^{k\to l})]=\beta_{kl}{\mathbf 1}(n_k\ge 1)&\mbox{jump of 
a task, for  $k\ne j\ne l$,}\\
{R}[(j,\unn) ;(j,\unn +\ue^{j\to l})]=\theta_{jl}{\mathbf 1}(n_j\ge 1)&\mbox{a task jumps from the loaded site, for $l\ne j$, }\\
{R}[(j,\unn) ;(j,\unn +\ue^{k\to j})]=\theta_{kj}e^\varphi{\mathbf 1}(n_k\ge 1)&\mbox{a task jumps
to the loaded site, for $k\ne j$,}\\
{R}[(j,\unn) ;(j^\prime, \unn )]=e^{-\varphi n_j}\tau_{jj^\prime}&\mbox{leap of a DC, from $j$ to $j^\prime\neq j$.} \end{array}\eeq

Pictorially, the DC is served or provides service affecting
task arrivals (both exogenous and intrinsic) at the loaded site, after which it moves to another 
site. Conversely, the DC departure from
the loaded site is influenced by the number of tasks accumulated in the queue. If $\vphi <0$ 
then the presence of the DC suppresses task arrivals, whereas  $\vphi >0$ means the opposite. 
Likewise, the
intensity of a DC leap increases with $n_j$ for $\vphi <0$ and decreases for $\vphi >0$.
\vskip .5 truecm

\section{An exact solution for a single DC} 

\subsection{The basic case}

In this sub-section we assume that rate collections $\rB$, $\Theta$ and $\rT$ satisfy 
\eqref{2.1}, \eqref{2.4} and are irreducible (i.e., matrices $\rB^u$, $\Theta^u$ and $\rT^u$ have
strictly positive entries for some positive integer $u$).

\begin{theorem}\label{t3.1} Assume that the sub-criticality condition holds true: 
$\;\diy\frac{\lam}{\mu},\;\frac{\lam e^\vphi}{\mu} <1$. 
Then the MP on $\Lam\times\Z_+^\Lam$ with generator $\bR$ is positive recurrent and
reversible (PRR). The stationary probability (SP) $\pi (j,\unn )$, $j\in\Lam$, $\unn =(n_i)\in\Z_+^\Lam$, 
of locating the DC at site $j$ and having $n_i$ tasks in sites of $\Lam$ is of the form 
\beq\label{eq:3.1}\pi (j,\unn )=\big(\#\,\Lam\big)^{-1}\left(1-\frac{\lam e^\vphi}{\mu}\right)
\left(1-\frac{\lam}{\mu}\right)^{\#\,\Lam\,-\,1}
\left(\frac{\lam}{\mu}\right)^{\sum\limits_{i\in\Lam}n_i} e^{\vphi n_j}.
\eeq 
\end{theorem}

\proof Probabilities $\pi (j,\unn )$ from \eqref{eq:3.1} satisfy the following DBEs:
\beq\label{eq:3.2}
\begin{array}{c}
\pi (j,\unn ){R}[(j,\unn) ;(j,\unn +\ue^k)]=\pi (j,\unn+\ue^k ){R}[(j,\unn+\ue^k) ;(j,\unn )],\;\;k\ne j,\\
\pi (j,\unn ){R}[(j,\unn) ;(j,\unn +\ue^j)]=\pi (j,\unn+\ue^j ){R}[(j,\unn+\ue^j) ;(j,\unn )],\\
\pi (j,\unn ){R}[(j,\unn) ;(j,\unn +\ue^{k\to l})]=\pi (j,\unn+\ue^{k\to l} ){R}[(j,\unn+\ue^{k\to l}) ;(j,\unn )],\;
k\neq j\neq l\neq k,\;n_k>0,\\
\pi (j,\unn ){R}[(j,\unn) ;(j,\unn +\ue^{j\to k})]=\pi (j,\unn+\ue^{j\to k}){R}[(j,\unn+\ue^{j\to k}) ;(j,\unn )],\; 
j\neq k,n_j>0,\\
\pi (j,\unn ){R}[(j,\unn) ;(j^\prime,\unn )]=\pi (j^\prime,\unn){R}[(j^\prime,\unn) ;(j,\unn )],\;j\neq j^\prime.
\end{array}\eeq
In fact, substituting \eqref{2.6}, omitting constant factors and canceling the common term $\diy\left(\frac{\lam}{\mu}\right)^{\sum_in_i}$ yields the identities: 
$$\begin{array}{c}
\diy e^{\varphi n_j}\lam =\left(\frac{\lam}{\mu}\right)e^{\varphi n_j}\mu,\;\;k\ne j,\;\;
\diy e^{\varphi n_j}\lam e^\varphi =\left(\frac{\lam}{\mu}\right) e^{\varphi (n_j+1)}\mu,\\
\diy e^{\varphi n_j}\beta_{kl} = e^{\varphi n_j}\beta_{lk},\;\;k\neq j\neq l\neq k,\;n_k>0,\\
\diy e^{\varphi n_j}\theta_{jk} = e^{\varphi (n_j-1)}\theta_{kj}e^\varphi ,\;\;j\neq k,n_j>0,\;\;
\diy e^{\varphi n_j}e^{-\varphi n_j}\tau_{jj^\prime} =e^{\varphi n_{j^\prime}}e^{-\varphi n_{j^\prime}}\tau_{j^\prime j},\;j\neq j^\prime.
\end{array}$$

Finally, under the irreducibility assumption, the process is PR. $\quad\Box$
\vskip .5 truecm

Note that the SP distribution $\pi$ in \eqref{eq:3.1} does not involve $\rB$, $\Theta$, $\rT$. The 
same pattern will be observed in the generalizations of the basic model below.

\subsection{A direct generalization}

In this sub-section, the rates $R\big[(j,\unn ),(j^\prime ,\unn^\prime )\big]$ are as follows (cf. \eqref{2.6}):
\beq\label{eq:3.3}\begin{array}{c}
{R}[(j,\unn) ;(j,\unn +\ue^k)]=\lam_k(n_k;j),\;k\neq j,\;\;
{R}[(j,\unn) ;(j,\unn +\ue^j)]=\lam_j(n_j;j)\gam_j(n_j),\\
{R}[(j,\unn) ;(j,\unn -\ue^k)]=\mu_k(n_k,k){\mathbf 1}(n_k\ge 1),\\
{R}[(j,\unn) ;(j,\unn +\ue^{k\to l})]=\beta_{kl}(n_k,n_l;j){\mathbf 1}(n_k\ge 1),\;k\neq j\neq l\neq k,\\
{R}[(j,\unn) ;(j,\unn +\ue^{j\to l})]=\theta_{jl}(n_j,n_l){\mathbf 1}(n_j\ge 1),\;l\neq j,\\
{R}[(j,\unn) ;(j,\unn +\ue^{k\to j})]=\theta_{kj}(n_k,n_j)\gam_j(n_j){\mathbf 1}(n_k\ge 1),\;k\neq j,\\
{R}[(j,\unn) ;(j^\prime,\unn )]=\big[\ogam_j(n_j)\big]^{-1}\tau_{jj^\prime}(\unn ),\;j\neq j^\prime . \end{array}\eeq

Let us comments on the form of these rates. As one can see, we take into account 
the numbers of tasks in related sites and the position of the DC. To start with, we deal 
in \eqref{eq:3.3} with values 
$\lam_i(n;j)$ and $\mu_i(n;j)$, $n\in\bbZ_+$, $i,j\in\Lam$. For $i\neq j$, values
$\lam_i(n;j)$, $\mu_i(n;j)$ give the intensities of exogenous arrival 
and exit of tasks, depending on $i$, the site location, $n (=n_i)$, the current number of tasks
at the site, and on $j$, the current loaded site. For $i=j$,  $\mu_j(n;j)$ yields an exit intensity 
from a loaded site while $\lam_j(n;j)$ represents a nominal arrival intensity which will be modified
via a gauge function $\gam_j(n)$. Examples are \def\rL{{\rm L}} \def\rU{{\rm U}}
$$
\lam_i(n;j)=\lam^\rU\1 (n< C),\;\;\;\mu_i(n,j)=\mu^\rU\min\,[n,K],\;j\neq i,$$
$$\lam_j(n;j)=\lam^\rL\1 (n< C),\;\;\;\;\mu_j(n;j)=\mu^\rL\min\,[n,K],$$ 
where $\lam^\rU$, $\lam^\rL$, $\mu^\rU$, $\mu^\rL$, $C$ and $K$ are positive constants
(which can be made varying the site). It means that the arrival
at a given (unloaded) site is blocked if the number of tasks reaches  $C$, and the pre-exit service is 
done by a $K$-server device, with intensities  depending on the site. 

We assume for simplicity that 
\beq\label{eq:3.4}\begin{array}{c}\sup\,\big[\lam_i(n;j):\;n\geq 0,\;i;j\in\Lam\big] <+\infty ,\;\;
\inf\,\big[\mu_i(n;j):\;n\geq 1,\;i;j\in\Lam\big] >0,\end{array}\eeq
and that $\lam_i(n;j)>0$ if $\lam_i(n+1;j)>0$.
Through the whole paper, we also set: 
\beq\label{eq:3.5}\olam_i(0;j)=\omu_i(0;j)=1,\;\;\;\olam_i(n;j)=\prod\limits_{0\leq m< n}\lam_i(m;j),
\;\omu_i(n;j)=\prod\limits_{1\leq m\leq n}\mu_i(m;j),\;n>0.\eeq

We also work with values $\gam_i(n)\geq 0$, $i\in\Lam$, $n\in\bbZ_+$, assuming 
that $\gam_i(0)=1$ and $\gam_i(n)>0$ if $\gam_i(n+1)>0$. These values are used to modify 
intensities of task arrival and jump at loaded site $j$. Viz., 
$$\gam_j(0)=1,\;\;\gam_j(n)=\diy\frac{e^{\vphi(j)}}{n}\;\hbox{ or }\;\gam_j(n)=ne^{\vphi(j)},\;n\geq 1,$$
where $\vphi (j)$ is a given real parameter depending upon $j$.  Next, we let, again throughout the paper,
\beq\label{eq:3.6}\ogam_i(0)=1,\;\;\;\ogam_i(n)=\prod\limits_{0\leq m<n}\gam_i(m),\;\;n\geq 1.
\eeq

Further, the intensities $\beta_{ik}$ depend on $j$ and $n_i$, $n_k$, and we write $\rB(\unn ;j)
=(\beta_{ik}(n_i,n_k;j))$. Consequently, we modify  symmetry assumptions in \eqref{2.1}: for $i\neq k$,
\beq\label{eq:3.7}\diy\frac{\lam_i(n_i-1;j)}{\mu_i(n_i;j)}\beta_{ik}(n_i,n_k;j)
=\frac{\lam_k(n_k;j)}{\mu_k(n_k+1;j)}\beta_{ki}(n_i+1,n_k-1;j),\eeq
and set $\beta_{ii}(n_i,n_i;j)=0$. Similar conditions are imposed on array $\Theta (\unn ;j)
=(\theta_{ik}(n_i,n_k;j))$: for $i\neq k$,
\beq\label{eq:3.8}
\diy\frac{\lam_i(n_i-1;j)}{\mu_i(n_i;j)}\theta_{ik}(n_i,n_k;j)=\frac{\lam_k(n_k;j)}{\mu_k(n_k+1;j)}
\theta_{ki}(n_k+1,n_i-1;j),\eeq
and $\theta_{ii}(n_i,n_i;j)=0$.

Finally, intensities $\tau_{ik}$ depend upon $\unn$ yielding an array $\rT (\unn) =
(\tau_{ik}(\unn),\;i,k\in\Lam )$. In this section we assume that $\tau_{ii}(\unn )=0$
and for $i\neq k$ and $\unn =(n_q,\,q\in\Lam)\in\Z_+^\Lam$,
\beq\label{eq:3.9}\tau_{ik}(\unn)\prod\limits_{q\in\Lam}\frac{\olam_q(n_q;i)}{\omu_q(n_q;i)}=\tau_{ki}
(\unn)\prod\limits_{q\in\Lam}\frac{\olam_q(n_q;k)}{\omu_q(n_q;k)}.\eeq
We continue referring to \eqref{eq:3.7} -- \eqref{eq:3.9} as symmetry conditions; one of primary 
future tasks should be
their replacement with less restrictive assumptions.  

In Theorem 3.2  we assume  conditions \eqref{eq:3.7} -- \eqref{eq:3.9} (these conditions 
will be recast in Section 4 in a more general situation) and suppose that
 arrays $\rB(\unn ;j)$, $\Theta (\unn ;j)$ and $\rT (\unn )$ have off-diagonal entries $>0$. 
 (We keep referring to this property as irreducibility.) The sub-criticality condition reads
\beq\label{eq:3.10}U_{ij}:=\sum\limits_{n\in\Z_+}\frac{\olam_i(n;j)}{\omu_i(n;j)}<+\infty,\;\;L_j:=\sum_{n\in\Z_+}
\frac{\olam_j(n;j)\ogam_j (n)}{\omu_j(n;j)} <\infty ,\;\;\forall\;\;i;j\in\Lam ,\;i\neq j.\eeq
 
\begin{theorem}\label{t3.2} Under \eqref{eq:3.10}, 
the MP on $\Lam\times\Z_+^\Lam$ with generator $\bR =\Big(
{R}[(j,\unn) ;(j^\prime,\unn^\prime)] \Big)$ as in \eqref{eq:3.3} is PRR. The SP $\pi (j,\unn )$,
of having the DC at site $j$ and $n_q$ tasks at $q\in\Lam$, is of the form 
\beq\label{eq:3.11}\pi (j,\unn )=\frac{1}{\Xi_\Lam}\,\prod\limits_{q\in\Lam}\,
\frac{\olam_q(n_q;j)}{\omu_q(n_q;j)}\,\ogam_j(n_j),\;j\in\Lam, \unn =(n_q)\in\Z_+^\Lam,\;\hbox{ with }\;\
\Xi_\Lam =\sum\limits_{j\in\Lam}L_j
\prod\limits_{i\in\Lam\setminus\{j\}}U_{ij}.\eeq 
\end{theorem}
Here $\;\Xi_\Lam\;$ is the partition function of the model.

\proof Probabilities $\pi (j,\unn )$ from \eqref{eq:3.11} and rates $R\big[(j,\unn ),(j^\prime ,\unn^\prime )\big]$ from
\eqref{eq:3.3} satisfy the DBEs \eqref{eq:3.2}.  In fact, after omitting the factor $\diy\frac{1}{ \Xi_\Lam}$ and canceling
common terms in $\;\diy\prod\limits_{i\in\Lam}\,\frac{\olam_i(n_i;j)}{\omu_i(n_i;j)}\;$, the DBEs are again 
verified with the help of symmetry conditions: for $j^\prime ,j,k,l\in\Lam$, with $k\neq j\neq l\neq k$, 
$j\neq\j^\prime$,
$$
\diy\ogam_j(n_j) \lam_k(n_k;j) =
\frac{\lam_k(n_k;j)\ogam_j(n_j)}{\mu_k(n_k+1;j)}
 \mu_k(n_k+1;j),$$
$$\diy\ogam_j(n_j) \lam_j(n_j;j)\gam_j(n_j)=
\frac{\lam_j(n_j;j)\ogam_j(n_j+1)}{\mu_j(n_j+1;j)}  \mu_j(n_j+1;j),$$ 
$$\begin{array}{l}
\diy\frac{\lam_k(n_k-1;j)}{\mu_k(n_k;j)}\ogam_j(n_j) \beta_{kl}(n_k,n_l;j) =
\frac{\lam_l(n_l;j)\ogam_j(n_j)}{\mu_l(n_l+1;j)} \beta_{lk}(n_l+1,n_k-1;j),\;
n_k\geq 1,\end{array}$$
$$\begin{array}{l}\diy\frac{\lam_j(n_j-1;j)\ogam_j(n_j)}{\mu_j(n_j;j)} \theta_{jk}(n_j,n_k;j) =\frac{\lam_k(n_k;j)\ogam_j(n_j-1)}{\mu_k(n_k+1;j)} \theta_{kj}(n_k+1,n_j-1;j)\gam_j(n_j-1) ,\;
n_j\geq 1,
\end{array}$$
$$\diy\ogam_j(n_j) \ogam_j(n_j)^{-1}\tau_{jj^\prime} (\unn )=\ogam_{j^\prime}(n_{j^\prime}) 
\ogam_{j^\prime}(n_{j^\prime})^{-1}\tau_{j^\prime j}(\unn ).$$ 
\vskip .3 truecm

As before, the process is PR under the irreducibility assumption.  
$\quad\Box$

\subsection{A closed-network version}

This version arises when we keep fixed the number of tasks in the network.  Correspondingly, we drop the 
two first lines in \eqref{eq:3.3}:
\beq\label{eq:3.12}\begin{array}{c}
{R}[(j,\unn) ;(j,\unn +\ue^{k\to l})]=\beta_{kl}(n_k,n_l;j){\mathbf 1}(n_k\ge 1),\;k\neq j\neq l,\\
{R}[(j,\unn) ;(j,\unn +\ue^{j\to l})]=\theta_{jl}(n_j,n_l){\mathbf 1}(n_j\ge 1),\\
{R}[(j,\unn) ;(j,\unn +\ue^{k\to j})]=\theta_{kj}(n_k,n_j)\gam_j(n_j){\mathbf 1}(n_k\ge 1),\\
{R}[(j,\unn) ;(j^\prime,\unn )]=\big[\ogam_j(n_j)\big]^{-1}\tau_{jj^\prime}(\unn ),\;j\neq j^\prime. \end{array}\eeq

In Theorem 3.3 we assume a modification of  condition \eqref{eq:3.9}
\beq\label{eq:3.13}\tau_{jj^\prime} (\unn )=\tau_{j^\prime j} (\unn ).\eeq

\medskip\begin{theorem}\label{t3.3} Fix $N$, the number of tasks in the network.
Given $\unn\in\Z_+^\Lam$, set: $|\unn |=\sum\limits_{s\in\Lam}n_s$.
The MP on $\left\{(j,\unn)\in\Lam\times\Z_+^\Lam :\;|\unn|=N\right\}$ with generator 
$\bR=\Big(R[(j,\unn) ;(j^\prime,\unn^\prime)] \Big)$ as in \eqref{eq:3.12} is PRR. The SPs 
$\;\pi (j,\unn )\;$ take the form 
$$\label{3.16}\pi (j,\unn )=\frac{\1(|\unn |=N)}{\Xi_{N,\Lam }}\;\ogam_j(n_j)\;\;\hbox{ where }\;\;
\Xi_{N,\Lam}=\sum\limits_{\substack{\unn=(n_s)\in\Z_+^\Lam :\\|\unn |=N}}
\;\sum\limits_{l\in\Lam}\ogam_l(n_l)\,.$$
\end{theorem} 

\proof Again, we use the DBEs verified with the help of the corresponding symmetry conditions:
$$\label{3.17}\begin{array}{c}
\pi (j,\unn ){R}[(j,\unn) ;(j,\unn +\ue^{i\to k})]=\pi (j,\unn+\ue^{i\to k}){R}[(j,\unn+\ue^{i\to k}) ;(j,\unn )],\;
i\neq j\neq k,n_i>0,\\
\pi (j,\unn ){R}[(j,\unn) ;(j,\unn +\ue^{j\to k})]=\pi (j,\unn+\ue^{j\to k}){R}[(j,\unn+\ue^{j\to k}) ;(j,\unn )],\;
j\neq k,n_j>0,\\
\pi (j,\unn ){R}[(j,\unn) ;(j^\prime,\unn )]=\pi (j^\prime,\unn){R}[(j^\prime,\unn) ;(j,\unn )],\;j\neq j^\prime.
\quad\Box\end{array}$$
 
\section{Simple exclusion in a Jackson-type environment} 

\subsection{A closed-open network}
\def\rE{{\rm E}}

The simple exclusion model was introduced in \cite{Spi1} (where the corresponding term  has 
been coined). The model was extensively studied thereafter: cf. \cite{Lig1}, \cite{Lig2}, \cite{Lig3}.  Here the state of the MP is a pair $(\uy,\unn)$ where $\uy =(y_s,s\in\Lam )\in\{0,1\}^\Lam$ and 
$\unn =(n_i,\,i\in\Lam)\in\Z_+^\Lam$. Let us set: $|\uy |=\sum\limits_{s\in\Lam}y_s$. We also write
$j\in\uy$ when $y_j=1$ and $j\not\in\uy$ when $y_j=0$.  In the case of a closed-open network, the 
sum $M=|\uy |$ remains constant. The rates (now denoted by ${R}[(\uy,\unn) ;(\uy^\prime,\unn^\prime)]$)
are specified as 
\beq\label{eq:4.1}\begin{array}{c}
{R}[(\uy,\unn) ;(\uy,\unn +\ue^p)]=\lam_p(n_p;\uy )[\gam_p(n_p)]^{y_p},\;
{R}[(\uy,\unn) ;(\uy,\unn -\ue^p)]=\mu_p(n_p;\uy ){\mathbf 1}(n_p\ge 1),\\
{R}[(\uy,\unn) ;(\uy,\unn +\ue^{k\to l})]=\beta_{kl}(n_k,n_l;\uy ),\;
{R}[(\uy,\unn) ;(\uy,\unn +\ue^{i\to j})]=\eps_{ij}(n_i,n_j;\uy ),\\
{R}[(\uy,\unn) ;(\uy,\unn +\ue^{i\to l})]=\theta_{il}(n_i,n_l;\uy ),\;
{R}[(\uy,\unn) ;(\uy,\unn +\ue^{k\to j})]=\theta_{kj}(n_k,n_j;\uy )\gam_j(n_j),\\
{R}[(\uy,\unn) ;(\uy +\ue^{j\to j^\prime},\unn )]=\big[\ogam_j(n_j)\big]^{-1}\tau_{jj^\prime}(\unn ;\uy), 
\end{array}\eeq
assuming (i) $i\neq j$, $i,j\in\uy$, $k\neq l$, $k,l,j^\prime\not\in\uy$, and (ii) $n_i,n_k\geq 1$.
Here we deal with intensities $\lam_\bullet (\,\cdot\,;\uy )$ and $\mu_\bullet (\,\cdot\,;\uy )$ depending on 
$\uy$. Such a generalization is extended to arrays 
$\rB (\unn;\uy)=(\beta_{kl}(n_k,n_l;\uy))$, $\Theta (\unn;\uy) =(\theta_{kl}(n_k,n_l;\uy ))$ and
$\rT (\unn, \uy )=(\tau_{jj^\prime}(\unn ;\uy ))$. We assume symmetry conditions similar to \eqref{eq:3.7},  
\eqref{eq:3.8} and \eqref{eq:3.9}: for $j\in\uy$, and $k\neq l$, $k,l,j^\prime\not\in\uy$,
\def\beac{\begin{array}{c}} \def\ena{\end{array}}
\beq\label{eq:4.2}
\beac\diy\frac{\lam_k(n_k-1;\uy)}{\mu_k(n_k;\uy)}\beta_{kl}(n_k,n_l;\uy)=\frac{\lam_l(n_l;\uy)}{\mu_l(n_l
+1;\uy)}\beta_{lk}(n_l+1,n_k-1;\uy),\;n_k\geq 1,\\
\diy\frac{\lam_j(n_j-1;\uy )}{\mu_j(n_j;\uy)}\theta_{jk}(n_j,n_k;\uy)=\frac{\lam_k(n_k;\uy)}{\mu_k(n_k+1;\uy)}
\theta_{kj}(n_k+1,n_j-1;\uy),\;\;n_j\geq 1,\ena\eeq
and -- see Eqn \eqref{eq:3.5} --
\beq\label{eq:4.3}
\diy \tau_{jj^\prime}(\unn ;\uy)\prod\limits_{q\in\Lam}\frac{\olam_q(n_q;\uy )}{\omu_q(n_q;\uy )}=\tau_{j^\prime j}
(\unn ;\uy+\ue^{j\to j^\prime})
\prod\limits_{q\in\Lam}\frac{\olam_q(n_q;\uy+\ue^{j\to j^\prime})}{\omu_q(n_q;\uy+\ue^{j\to j^\prime})}. \eeq
 \def\beal{\begin{array}{l}} \def\bear{\begin{array}{r}}

We also have a new collection of jump rates $\rE (\unn;\uy)=(\eps_{ij}(n_i,n_j;\uy))$ satisfying the 
symmetry property: for $i\neq j$, $i,j\in\uy$ and $n_i\geq 1$,
\beq\label{eq:4.4}\beal
\diy\frac{\lam_i(n_i-1;\uy )\gam_i(n_i-1)}{\mu_i(n_i;\uy)}\eps_{ij}(n_i,n_j;\uy)
=\frac{\lam_j(n_j;\uy)\gam_j(n_j)}{\mu_j(n_j+1;\uy)}\eps_{ji}(n_j+1,n_i-1;\uy).
\ena\eeq

Until the end of Section 4 we work with
irreducible collections $\rB(\unn ;\uy )$, $\Theta (\unn ;\uy )$, $\rE(\unn ;\uy )$ and $\rT (\unn ;\uy )$. 
Furthermore, assumptions \eqref{eq:4.2} -- \eqref{eq:4.4} are adopted in sub-Sections 4.1 and 4.2.  
The interpretation is that we have a $1$-$0$ configuration $\uy$ of loaded sites occupied by DCs,
with a total number of DCs $M$; each of them influences task arrivals and task jumps as
indicated, independently for different sites. In addition, each DC can jump from a loaded to an unloaded
site, again independently.  

With $\ogam_j(n)$ as in \eqref{eq:3.6}, the SPs $\pi (\uy,\unn )$ read
\beq\label{eq:4.5}\pi (\uy,\unn )=\frac{{\mathbf 1}(|\uy |=M)}{\Xi_{\Lam,M}}\prod\limits_{q\in\Lam}
\frac{\olam_q(n_q;\uy )}{\omu_q(n_q;\uy )}\prod\limits_{j\in\Lam:\;y_j=1}\ogam_j(n_j)
\eeq 
with partition function $\diy\Xi_{\Lam ,M}=\sum\limits_{\substack{\uy=(y_s)\in\{0,1\}^\Lam :\\
|\uy |=M}}\;\;\prod\limits_{r:\;y_r=1}L_r(\uy )
\prod\limits_{l:\;y_l=0}U_l(\uy )\,,\;\;M\leq\#\,\Lam$,
and
\beq\label{eq:4.6}
U_l(\uy )=\sum\limits_{n\in\Z_+}\frac{\olam_l(n;\uy )}{\omu_l(n;\uy )},\;\;\;L_r(\uy )=\sum\limits_{n\in\Z_+}\frac{\olam_r(n;\uy )\ogam_r(n)}{\omu_r(n;\uy )}.
\eeq 
The sub-criticality conditions emerging from \eqref {eq:4.6} are: $\forall$ $\uy\in\{0,1\}^\Lam$ with $|\uy|=M$,
\beq\label{eq:4.7}\begin{array}{lr}
\diy U_l(\uy )<+\infty ,\;L_j(\uy )<+\infty ,&\hbox{$\forall$ $l,j\in\Lam$ with $y_l=0$ and $y_j=1$.}\end{array}
\eeq 

\begin{theorem}\label{t4.1} The MP with generator $\bR=\Big(
{R}[(\uy,\unn) ;(\uy^\prime,\unn^\prime)] \Big)$ as in {\rm{\eqref{eq:4.1}}} on state space \\
$\Big\{(\uy,\unn)\in\{0,1\}^\Lam
\times\Z_+^\Lam :\;|\uy|=M\Big\}$  is PRR. The SPs $\pi (\uy,\unn )$ are given by \eqref{eq:4.5}. 
\end{theorem}

\proof As before, the proof is based on DBEs. These are now as follows: for $j,j^\prime ,k,l\in\Lam$,
with $k\neq l$, 
\beq\label{eq:4.8}
\begin{array}{c}
\pi (\uy,\unn ){R}[(\uy,\unn) ;(\uy,\unn +\ue^k)]=\pi (\uy,\unn+\ue^k ){R}[(\uy,\unn+\ue^k) ;(\uy,\unn )],
\;\;k\not\in\uy,\\
\pi (\uy,\unn ){R}[(\uy,\unn) ;(\uy,\unn +\ue^j)]=\pi (\uy,\unn+\ue^j ){R}[(\uy,\unn+\ue^j) ;(\uy,\unn )\,],
\;\;j\in\uy,\\
\pi (\uy,\unn ){R}[(\uy,\unn) ;(\uy,\unn +\ue^{k\to l})]=\pi (\uy,\unn+\ue^{k\to l}){R}
[(\uy,\unn +\ue^{k\to l}) ;(\uy,\unn )],\;
y_k=y_l,\;n_k\geq 1,\\
\pi (\uy,\unn ){R}[(\uy,\unn) ;(\uy,\unn +\ue^{j\to k})]=\pi (\uy,\unn +\ue^{j\to k}){R}
[(\uy,\unn +\ue^{j\to k}) ;(\uy,\unn )],\;
j\in\uy,\;k\not\in\uy,\;n_j\geq 1,\\
\pi (\uy,\unn ){R}[(\uy,\unn) ;(\uy+\ue^{j\to j^\prime},\unn )]=\pi (\uy+\ue^{j\to j^\prime},\unn)
{R}[(\uy +\ue^{j\to j^\prime},\unn);(\uy,\unn )],\;
j\in\uy,\;j^\prime\not\in\uy.
\end{array}\eeq
The verification is still direct. For definiteness, we show the equation emerging in the third line of 
\eqref{eq:4.8}, when $k,l\in\uy$ (other cases have been effectively considered earlier). Upon 
omitting the factor $\diy\frac{1}{ \Xi_\Lam}$ and canceling 
common terms in the products $\diy\prod\limits_{j\in\Lam:\;y_j=1}\ogam_j(n_j)$ and $\diy
\prod\limits_{q\in\Lam}\,\frac{\olam_q(n_q;j)}{\omu_q(n_q;j)}$, this equation becomes
\eqref{eq:4.4}. $\quad\Box$

\subsection{An open-open network}

Now the rates from \eqref{eq:4.1} are complemented with
\beq\label{eq:4.9}\begin{array}{ll}
{R}[(\uy,\unn) ;(\uy +\ue^k,\unn )]=\xi_k\ogam_k(n_k)\1(y_k=0),&\hbox{arrival of a DC at site $k$,}\\
{R}[(\uy,\unn) ;(\uy -\ue^i,\unn )]=\eta_i\1(y_i=1),&\hbox{exit of a DC from site $i$.}\end{array}
\eeq
Here $\xi_k>0$ and $\eta_i>0$ do not depend on $\unn$ or $\uy$. (This assumption can be 
generalized but certain independence should be maintained.) Further, in this sub-section we assume 
Eqn \eqref{eq:4.2}. In addition, we assume here that, $\forall$ $\unn\in\bbZ_+^\Lam$, 
\beq\label{eq:4.10}\beal
\hbox{the product }\;
V(\unn )=\diy\prod\limits_{q\in\Lam}
\frac{\olam_q(n_q;\uy )}{\omu_q(n_q;\uy )}\; 
\hbox{ does not depend on configuration $\uy\in\{0,1\}^\Lam$.}
\ena\eeq
Accordingly, Eqn \eqref{eq:4.3} has to be replaced with 
\beq\label{eq:4.11} \tau_{jj^\prime}(\unn ;\uy)\frac{\xi_j}{\eta_j}
=\tau_{j^\prime j}(\unn ;\uy+\ue^{j\to j^\prime})\frac{\xi_{j^\prime}}{\eta_{j^\prime}},\;\;j\in\uy,\;j^\prime\not\in\uy.
\eeq

The SPs $\pi (\uy,\unn )$ become
\beq\label{eq:4.12}\pi (\uy,\unn )=\frac{V(\unn )}{\Xi_\Lam}
\prod\limits_{j\in\Lam:\;y_j=1}\frac{\xi_j\ogam_j(n_j)}{\eta_j}.
\eeq 
Here $\diy\Xi_\Lam = \sum\limits_{\uy\in\{0,1\}^\Lam} \sum\limits_{\unn\in\bbZ_+^\Lam}
V(\unn )\prod\limits_{j\in\Lam:\;y_j=1}\frac{\xi_j\ogam_j(n_j)}{\eta_j}
= \sum\limits_{\uy\in\{0,1\}^\Lam}\;\prod\limits_{r:y_r=1}\frac{\xi_rL_r(\uy)}{\eta_r}
\prod\limits_{l\in\Lam :\;y_l=0}U_l(\uy )$, 
and $L_r(\uy)$ and $U_l(\uy)$ are as in \eqref{eq:4.6}. The sub-criticality condition reads $ \Xi_\Lam<\infty$, or 
\beq\label{eq:4.13}\begin{array}{r}
\diy U_l(\uy )<+\infty,\;L_r(\uy )<+\infty ,\;\hbox{ $\forall$ $\uy\in\{0,1\}^\Lam$}\;
\hbox{and $l,r\in\Lam$ with $y_l=0$ and $y_r=1$,}\end{array}
\eeq 
and can be treated as a slight modification of \eqref{eq:4.7}. 
\vskip .5 truecm

\begin{theorem}\label{t4.2} Under \eqref{eq:4.13}, 
the MP on $\{0,1\}^\Lam\times\Z_+^\Lam$ with generator $\bR=\Big(
{R}[(\uy,\unn) ;(\uy^\prime,\unn^\prime)] \Big)$ as in {\rm{\eqref{eq:4.1}, \eqref{eq:4.8}}} is PRR. The SPs 
$\;\pi (\uy,\unn )\;$ are given by \eqref{eq:4.12}. 
\end{theorem}

\proof As before, one checks the DBEs \eqref{eq:4.8} completed with
\beq\label{eq:4.14}  \pi (\uy,\unn ){R}[(\uy,\unn) ;(\uy +\ue^k,\unn )]=
 \pi (\uy +\ue^k,\unn ){R}[(\uy +\ue^k,\unn );(\uy,\unn)],\;\;k\not\in\uy .
\eeq
In view of \eqref{eq:4.10}, the latter holds true. $\qquad\Box$

\subsection{A closed-closed network}

Here we keep $|\uy |$ and $|\unn |$ fixed: $|\uy |=M$, $|\unn |=N$. The rates are
as in Eqn \eqref{eq:4.1}, with the top three lines discarded. The following conditions are assumed: for $i\neq j$, $i,j\in\uy$ and $k\neq l$, $k,l,j^\prime\not\in\uy$,
\beq\label{eq:4.15}\beac
\beta_{kl}(n_k,n_l;\uy)=\beta_{lk}(n_l+1,n_k-1;\uy),\;
\theta_{jk}(n_j,n_k;\uy)=\theta_{kj}(n_k+1,n_j-1;\uy),\;n_k,n_j\geq 1,\\
\eps_{ij}(n_i,n_j;\uy) \gamma_i(n_i-1)=\gamma_j(n_j) \eps_{ji}(n_j+1,n_i-1;\uy),\;\;n_i\geq 1, 
\ena\eeq
and
\beq\label{eq:4.16}
 \tau_{jj^\prime}(\unn ;\uy)=\tau_{j^\prime j}(\unn ;\uy-\ue^j+\ue^{j^\prime}),
\eeq
which can be viewed as a modification of  \eqref{eq:4.2} -- \eqref{eq:4.4} and \eqref{eq:3.13}.

The SP distribution resembles \eqref{eq:3.11}:
\beq\label{eq:4.17}  \pi (\uy,\unn )=\frac{\1(|\uy |=M,\;|\unn |=N)}{\Xi_{N,\Lam ,M}}
\prod\limits_{j\in\Lam}\big[\ogam_j(n_j)\big]^{y_j},\;\;
\Xi_{N,\Lam ,M}
=\;\sum\limits_{\substack{\uy=(y_s)\in\{0,1\}^\Lam,\\ \unn=(n_s)\in\Z_+^\Lam :\\ |\uy|=M, |\unn |=N}}
\;\;\;\prod\limits_{l\in\Lam}\big[\ogam_l(n_l)\big]^{y_l}\,.\eeq 

\begin{theorem}\label{t4.3} 
Given integer $M,N\geq 1$, the MP on $\Big\{(\uy,\unn)\in\{0,1\}^\Lam
\times\Z_+^\Lam :\;|\uy|=M,|\unn|=N\Big\}$ with generator $\bR=\Big(
{R}[(\uy,\unn) ;(\uy^\prime,\unn^\prime)] \Big)$ as specified in this sub-section is PRR. The SPs
$\pi (\uy,\unn )$ are given by {\rm{\eqref{eq:4.17}}}. 
\end{theorem}

\proof Again by means of suitable DBEs \eqref{eq:4.8}, with the help of \eqref{eq:4.15}. $\quad\Box$

\subsection{An open-closed network}

In this version of the model, $N$, the number of tasks, is fixed, but the number of DC`s varies due to arrivals and exits. A part of transition rates $R[(\uy,\unn );(\uy^\prime,\unn^\prime)]$ comes from \eqref{eq:4.1}:
\beq\label{eq:4.18}\begin{array}{c}
{R}[(\uy,\unn) ;(\uy,\unn +\ue^{k\to l})]=\beta_{kl}(n_k,n_l;\uy ){\mathbf 1}(n_k\ge 1),\;k\neq l,\;k,l\not\in\uy,\\
{R}[(\uy,\unn) ;(\uy,\unn +\ue^{i\to j})]=\eps_{ij}(n_i,n_j;\uy ){\mathbf 1}(n_i\ge 1),\;i\neq j,\;i,j\in\uy,\\
{R}[(\uy,\unn) ;(\uy,\unn +\ue^{j\to l})]=\theta_{jl}(n_j,n_l;\uy ){\mathbf 1}(n_j\ge 1),\;j\in\uy,l\not\in\uy,\\
{R}[(\uy,\unn) ;(\uy,\unn +\ue^{k\to j})]=\theta_{kj}(n_k,n_j;\uy )\gam_j(n_j){\mathbf 1}(n_k\ge 1),\;j\in\uy,
k\not\in\uy,\\
{R}[(\uy,\unn) ;(\uy +\ue^{j\to j^\prime},\unn )]=\big[\ogam_j(n_j)\big]^{-1}\tau_{jj^\prime}(\unn ;\uy),\;j\in\uy,
j^\prime\not\in\uy . \end{array}\eeq
In addition, we use the rates from \eqref{eq:4.9}. The DBEs read: for $i,j,j^\prime,k,l\in\Lam$ with $k\neq l$,
$j\in\uy$,
\beq\label{eq:4.19}\beac\pi (\uy,\unn ){R}[(\uy,\unn) ;(\uy,\unn+\ue^{k\to l})]=\pi (\uy,\unn +\ue^{k\to l}){R}
[(\uy,\unn+\ue^{k\to l}) ;(\uy,\unn )],\;
y_k=y_l,\;n_k\geq 1,\\
\pi (\uy,\unn ){R}[(\uy,\unn) ;(\uy,\unn +\ue^{j\to k})]=\pi (\uy,\unn +\ue^{j\to k}){R}
[(\uy,\unn +\ue^{j\to k}) ;(\uy,\unn )],\; k\not\in\uy,\;n_j\geq 1,\\
\pi (\uy,\unn ){R}[(\uy,\unn) ;(\uy +\ue^{j\to j^\prime},\unn )]=\pi (\uy +\ue^{j\to j^\prime},\unn)
{R}[(\uy +\ue^{j\to j^\prime},\unn);(\uy,\unn )],\;j^\prime\not\in\uy,\\
\pi (\uy,\unn ){R}[(\uy,\unn) ;(\uy +\ue^k,\unn )]=\pi (\uy+\ue^k,\unn ){R}[(\uy +\ue^k,\unn );(\uy,\unn)],
\;\;k\not\in\uy .\end{array}\eeq

We now assume conditions \eqref{eq:4.11} and \eqref{eq:4.15}. The SPs and sub-criticality condition read
\beq\label{eq:4.20}
\pi (\uy,\unn )=\frac{\1(|\unn |=N)}{\Xi_{N,\Lam}}
\prod\limits_{q\in\Lam}\left[\frac{\xi_q}{\eta_q}\ogam_q(n_q)\right]^{y_q},\;\;\Xi_{N,\Lam}
=\;\sum\limits_{\substack{\uy=(y_s)\in\{0,1\}^\Lam,\\ \unn=(n_s)\in\Z_+^\Lam :\,|\unn |=N}}
\;\;\prod\limits_{q\in\Lam}\left[\frac{\xi_q}{\eta_l}\ogam_q(n_q)\right]^{y_q}<\infty\,.\eeq 
\vskip .5 truecm

\begin{theorem}\label{t4.4} 
Fix $N\in\bbZ_+$. If $\Xi_{N,\Lam} <\infty$, 
the MP with generator 
$\bR =\Big(R[(\uy;\unn); (\uy^\prime ;\unn^\prime)]\Big)$ as above is PRR on
state space $\Big\{(\uy,\unn )\in\{0,1\}^\Lam\times\Z_+^\Lam :\;|\unn |=N\Big\}$ . The 
SPs are given by \eqref{eq:4.20}.
\end{theorem}

\proof Still the DBEs, now from Eqn \eqref{eq:4.19}. $\quad\Box$

\section{A zero-range system in a Jackson-type environment} 

A zero-range modification arises when we allow the DCs to accumulate in sites
$i\in\Lam$. Here $\uy =(y_s,\,s\in\Lam)\in\Z_+^\Lam$; we again set  
$|\uy |=\sum\limits_{s\in\Lam}y_s$ and write $j\in\uy$ when $y_j\geq 1$ and $l\not\in\uy$
when $y_l=0$. 

\subsection{A closed-open network}

In this sub-section, $M:=|\uy |$ is a conserved quantity.
The rates are similar to \eqref{eq:4.1}:  for $i,j,k,l,p\in\Lam$,
\beq\label{5.1}\begin{array}{c}
{R}[(\uy,\unn) ;(\uy,\unn +\ue^p)]=\lam_p(n_p;\uy )\big[\gam_p(n_p)\big]^{y_p},\;\;
{R}[(\uy,\unn) ;(\uy,\unn -\ue^p)]=\mu_p(n_p;\uy ){\mathbf 1}(n_p\ge 1),\\
{R}[(\uy,\unn) ;(\uy,\unn +\ue^{k\to l})]=\beta_{kl}(n_k,n_l;\uy ){\mathbf 1}(n_k\ge 1),\;k,l\not\in\uy,\;k\neq l\\
{R}[(\uy,\unn) ;(\uy,\unn +\ue^{i\to j})]=\eps_{ij}(n_i,n_j;\uy ){\mathbf 1}(n_i\ge 1),\;i,j\in\uy,\;i\neq j\\
{R}[(\uy,\unn) ;(\uy,\unn +\ue^{j\to l})]=\theta_{jl}(n_j,n_l;\uy ){\mathbf 1}(n_j\ge 1),\;j\in\uy,\;l\not\in\uy ,\\
{R}[(\uy,\unn) ;(\uy,\unn +\ue^{k\to j})]=\theta_{kj}(n_k,n_j;\uy )\big[\gam_j(n_j)\big]^{y_j}
{\mathbf 1}(n_k\ge 1),\;j\in\uy,\;k\not\in\uy,\\
{R}[(\uy,\unn) ;(\uy +\ue^{j\to j^\prime},\unn )]=\big[\ogam_j(n_j)\big]^{-y_j}\tau_{jj^\prime}(\unn;\uy )\big[\ogam_{j^\prime}(n_{j^\prime})\big]^{-y_{j^\prime}},\;\;j\in\uy ,jö\prime\neq j. \end{array}\eeq
Until the end of Section 5, we suppose irreducibility of collections $\rB(\unn ;\uy )$, $\Theta (\unn ;\uy )$, $\rE(\unn;\uy )$ and $\rT (\unn ;\uy )$. 
In sub-Sections 5.1 and 5.2 intensities $\eps_{ij}(n_i,n_j;\uy)$ are supposed to obey: $\forall$ $i\neq j$, $i,j\in\uy$,
\beq\label{5.2}\beal
\diy\frac{\lam_i(n_i-1;\uy )\big[\gam_i(n_i-1)\big]^{y_i}}{\mu_i(n_i;\uy)}\eps_{il}(n_i,n_j;\uy)
=\frac{\lam_j(n_j;\uy)\big[\gam_j(n_j)\big]^{y_j}}{\mu_j(n_j+1;\uy)}\eps_{ji}(n_j+1,n_i-1;\uy),\;n_i\geq 1,\ena\eeq
which replaces \eqref{eq:4.4}.  Conditions \eqref{eq:4.2}--\eqref{eq:4.3} remain in place in
both sub-Sections 5.1 and 5.2.

The present model gives rise to SPs 
\beq\label{5.3}
\pi (\uy,\unn )=\frac{\1 (|\uy |=M)}{\Xi_{\Lam ,M}}\prod\limits_{p\in\Lam}
\frac{\olam_p(n_p;\uy )}{\omu_p(n_p;\uy )}
\prod\limits_{q\in\Lam}\big[\ogam_q(n_q)\big]^{y_q},\;
\diy\Xi_{\Lam ,M}=
\sum\limits_{\substack{\uy=(y_s)\in\Z_+^\Lam :\\
|\uy |=M}}\prod\limits_{r:\;y_r\geq 1}C_r(\uy )
\prod\limits_{l:\;y_l=0}U_l(\uy),\eeq 
and 
\beq\label{5.4}
U_l(\uy )=\sum\limits_{n\in\Z_+}\frac{\olam_l(n_l;\uy )}{\omu_l(n;\uy )},\;C_r(\uy )=\sum\limits_{n\geq 0}\frac{\olam_r(n;\uy )}{\omu_r(n;\uy )}\ogam_r(n)^{y_r}.
\eeq 
The sub-criticality conditions read: $\forall$ $\uy\in\Z_+^\Lam$ with $|\uy|=M$, 
\beq\label{5.5}\begin{array}{r}
\diy U_l(\uy )<+\infty ,\;C_r(\uy )<+\infty ,\;
\hbox{ $\forall$ $l,r\in\Lam$ with $y_l=0$ and $y_r\geq 1$.}\end{array}
\eeq 

\begin{theorem}\label{t5.1} Under conditions {\rm{\eqref{5.5}}},  
the MP on $\Big\{(\uy,\unn)\in\Z_+^\Lam
\times\Z_+^\Lam :\;|\uy|=M\Big\}$  with generator $\bR=\Big(
{R}[(\uy,\unn) ;(\uy^\prime,\unn^\prime)] \Big)$ as in {\rm{\eqref{5.1}}} is PRR. The SPs
$\pi (\uy,\unn )$ are given by {\rm{\eqref{5.3}, \eqref{5.4}}}. 
\end{theorem}

\proof Still the DBEs, now for rates \eqref{5.1}. The DBEs are as in \eqref{eq:4.8}, and after
cancellations are reduced to  \eqref{eq:4.2}, \eqref{eq:4.3} and \eqref{5.2}. $\quad\Box$

\subsection{An open-open network}

In this model we allow both the tasks and DCs to come and leave. Correspondingly, the rates 
\eqref{5.1} are complemented in a manner similar to \eqref{eq:4.9}:
\beq\label{5.6}\begin{array}{l}
{R}[(\uy,\unn) ;(\uy +\ue^p,\unn )]=\xi_p\ogam_p(n_p),\;\;
{R}[(\uy,\unn) ;(\uy -\ue^p,\unn )]=\eta_p\1(y_p\geq 1).\end{array}
\eeq

As in sub-Section 4.2, we assume \eqref{eq:4.2} and replace \eqref{eq:4.3} by \eqref{eq:4.10}.
We also assume \eqref{eq:4.11}. (In \eqref{eq:4.10}, $\uy\in\{0,1\}^\Lam$ is replaced with
$\uy\in\bbZ_+^\Lam$, and in \eqref{eq:4.11} the restriction $j^\prime\not\in\uy$ removed.)
As in sub-Section 5.1, condition \eqref{5.2} is also in place.

Assuming $U_l(\uy)$ and $C_l(\uy )$ as in \eqref{5.4}, the SPs and sub-criticality condition now read
\beq\label{5.7}
\pi (\uy,\unn )=\Xi_{\Lam}^{-1}\prod\limits_{p\in\Lam}
\frac{\olam_p(n_p;\uy )}{\omu_p(n_p;\uy )}
\prod\limits_{q\in\Lam}\left[\frac{\xi_q}{\eta_q}\ogam_q(n_q)\right]^{y_q}, \;\;
\diy\Xi_{\Lam}=
\sum\limits_{\uy=(y_s)\in\Z_+^\Lam}\;\prod\limits_{r:\;y_r\geq 1}C_r(\uy )
\prod\limits_{l:\;y_l=0}U_l(\uy)<\infty.\eeq

\begin{theorem}\label{t5.2} 
If $\;\Xi_\Lam <+\infty$,  
the MP on $\Z_+^\Lam\times\Z_+^\Lam $ with generator $\bR=\Big(
{R}[(\uy,\unn) ;(\uy^\prime,\unn^\prime)] \Big)$ as in {\rm{\eqref{5.1}}}, {\rm{\eqref{5.6}}} is PRR. The 
SPs $\;\pi (\uy,\unn )$ are given by {\rm{\eqref{5.7}}}. \end{theorem}

\proof The DBEs again. The added equations \eqref{5.6} are treated similarly to \eqref{eq:4.14}. \quad$\Box$
\vskip .5 truecm

\subsection{A closed-closed network}

Let us now suppose that both $|\uy |$ and $|\unn |$ are fixed: $|\uy |=M$ and $|\unn |=N$. The 
rates are as in Eqn \eqref{5.1}, with top two lines discarded. In this sub-section, we assume conditions \eqref{eq:4.15} and \eqref{eq:4.16}, with specification $\uy\in\{0,1\}$ replaced by $\uy\in\bbZ_+^\Lam$
and condition $j^\prime\not\in\uy$ removed.  The only exception is that
the bottom line in \eqref{eq:4.15} is now replaced with 
\beq\label{5.8}\eps_{ij}(n_i,n_j;\uy) [\gamma_i(n_i-1)]^{y_i}
=[\gamma_j(n_j)]^{y_j}\eps_{ji}(n_j+1,n_i-1;\uy),\;\; i\neq j,\;i,j\in\uy,\;n_i\geq 1. \eeq

The SP distribution mimicks \eqref{eq:4.17}:  
\beq\label{5.9}  \pi (\uy,\unn )=\frac{\1(|\uy |=M,|\unn |=N)}{\Xi_{N,\Lam ,M}}\prod\limits_{j\in\Lam}
\big[\ogam_j(n_j)\big]^{y_j},\;\;
\Xi_{N,\Lam ,M}
=\sum\limits_{\substack{\uy=(y_s),\unn=(n_s)\in\Z_+^\Lam :\\ |\uy|=M,\;|\unn |=N}}
\;\prod\limits_{l\in\Lam}\big[\ogam_l(n_l)\big]^{y_l}.\eeq 

The above DBEs and symmetry conditions (including \eqref{5.8}) lead to Theorem \ref{t5.3}:

\begin{theorem}\label{t5.3} 
The MP on state space $\Big\{(\uy,\unn)\in\Z_+^\Lam
\times\Z_+^\Lam :\;|\uy|=M,\;|\unn|=N\Big\}$ with generator\\ $\bR=\Big(
{R}[(\uy,\unn) ;(\uy^\prime,\unn^\prime)] \Big)$ as in {\rm{\eqref{5.1}}} is PRR. The SPs
$\pi (\uy,\unn )$ are given by {\rm{\eqref{5.9}}}. 
\end{theorem}

\subsection{An open-closed network}

Here -- as in sub-Section 4.4 -- we only fix $N$. The rates follow \eqref{eq:4.18} and consists of
\beq\label{5.10}\begin{array}{c}
{R}[(\uy,\unn) ;(\uy,\unn +\ue^{k\to l})]=\beta_{kl}(n_k,n_l;\uy ){\mathbf 1}(n_k\ge 1),\;k\neq l,\;k,l\not\in\uy,\\
{R}[(\uy,\unn) ;(\uy,\unn +\ue^{i\to j})]=\eps_{ij}(n_i,n_j;\uy ){\mathbf 1}(n_i\ge 1),\;i\neq j,\;i,j\in\uy,\\
{R}[(\uy,\unn) ;(\uy,\unn +\ue^{j\to l})]=\theta_{jl}(n_j,n_l;\uy ){\mathbf 1}(n_j\ge 1),\;j\in\uy,\;l\not\in\uy ,\\
{R}[(\uy,\unn) ;(\uy,\unn +\ue^{k\to j})]=\theta_{kj}(n_k,n_j;\uy )\big[\gam_j(n_j)\big]^{y_j}
{\mathbf 1}(n_k\ge 1),\;j\in\uy,\;k\not\in\uy,\\
{R}[(\uy,\unn) ;(\uy +\ue^{j\to j^\prime},\unn )]=\big[\ogam_j(n_j)\big]^{-y_j}\tau_{jj^\prime}(\unn;\uy )
\big[\ogam_{j^\prime}(n_{j^\prime})\big]^{-y_{j^\prime}},\;j\neq j^\prime ,\;j\in\uy ,\end{array}\eeq
plus rates from Eqns \eqref{5.6}. 
As in sub-Section 4.4, we now assume conditions \eqref{eq:4.11} and \eqref{eq:4.15}, modified like above (including \eqref{5.9}). The SPs and sub-criticality condition read
\beq\label{5.11}
\pi (\uy ,\unn )\;=\;\frac{{\mathbf 1}(|\unn |=N)}{\Xi_{N,\Lam}}\;\prod\limits_{l\in\Lam}\left[\frac{\xi_l}{\eta_l}
\ogam_l(n_l)\right]^{y_l},\;\;\Xi_{N,\Lam}
=\;\sum\limits_{\substack{\uy=(y_s) \unn=(n_s)\in\Z_+^\Lam :\\ |\unn |=N}}
\;\;\prod\limits_{l\in\Lam}\left[\frac{\xi_l}{\eta_l}\ogam_l(n_l)\right]^{y_l}<\infty \,.\eeq 
The DBEs in this case yield

\begin{theorem}\label{t5.4} Fix $N\in\bbZ_+$ and consider  the MP on $\Big\{(\uy,\unn)\in\Z_+^\Lam
\times\Z_+^\Lam :\;|\unn|=N\Big\}$ with generator $\bR=\Big(
{R}[(\uy,\unn) ;(\uy^\prime,\unn^\prime)] \Big)$ as in {\rm{\eqref{5.10}}}. If $\;\Xi_{N,\Lam}<\infty$, Assuming {\rm{\eqref{5.11}}}, 
it is PRR. The SPs are given by \eqref{5.11}.
\end{theorem}

\vskip .2 truecm

{\bf Acknowledgements}. We thank the referee and V. Belitsky for corrections/suggestions and express gratitude: MG to CAPES, CNPq; YS to USP, ICMC, Math Dept PSU; AY to FAPESP, CNPq, Dept of Pharmacy OSU, for support/hospitality. The research of E. Pechersky was carried out
at the IITP RAS at the expense of the Russian Foundation for Sciences (project
14-50-00150).

\vskip .5 truecm

M Gannon\quad mark@ime.usp.br

E Pechersky\quad pech@iitp.ru; epechersky@gmail.com

Y Suhov \quad yms@statslab.cam.ac.uk; ims14@psu.edu\quad (communicating author)

A Yambartsev\quad 	yambar@ime.usp.br; yambar@gmail.com
\end{document}